\documentclass[article,10pt,letterpaper,twoside]{article}

\usepackage[letterpaper,left=1in,right=1in,top=1.5in,bottom=1.5in]{geometry}

\title{The topological period--index conjecture}
\author{Benjamin Antieau and Ben Williams}
\date{\today}

\usepackage[pdfstartview=FitH,
 pdfauthor={Antieau and Williams},
 pdftitle={The topological period--index conjecture},
 colorlinks,
 linkcolor=reference,
 citecolor=citation,
 urlcolor=e-mail,
 backref]{hyperref}

\usepackage{amsmath}
\usepackage{amscd}
\usepackage{amsbsy}
\usepackage{amssymb}
\usepackage{verbatim}
\usepackage{eucal}
\usepackage{microtype}
\usepackage{hyperref}
\usepackage{mathrsfs}
\usepackage{amsthm}
\usepackage{stmaryrd}
\usepackage{xy}
\input xy
\xyoption{all}
\usepackage{tikz}
\usetikzlibrary{matrix,arrows}
\usepackage[abbrev,lite,alphabetic]{amsrefs}
\usepackage{dsfont}
\usepackage[T1]{fontenc}
\usepackage{mathbbol}

\usepackage{fancyhdr}
\usepackage{makeidx}
\makeindex

\usepackage{color}
\definecolor{todo}{rgb}{1,0,0}
\definecolor{conditional}{rgb}{0,1,0}
\definecolor{e-mail}{rgb}{0,.40,.80}
\definecolor{reference}{rgb}{.20,.60,.22}
\definecolor{mrnumber}{rgb}{.80,.40,0}
\definecolor{citation}{rgb}{0,.40,.80}

\renewcommand{\bf}{\bfseries}

\let\oldmarginpar\marginpar
\renewcommand\marginpar[1]{\-\oldmarginpar[\raggedleft\footnotesize #1]%
{\raggedright\footnotesize #1}}

\newcommand{\Ascr}{\mathcal{A}}

\newcommand{\E}{\mathrm{E}}

\renewcommand{\H}{\mathrm{H}}

\newcommand{\K}{\mathrm{K}}

\newcommand{\CC}{\mathds{C}}

\newcommand{\ZZ}{\mathds{Z}}

\newcommand{\tors}{\mathrm{tors}}

\renewcommand{\geq}{\geqslant}
\renewcommand{\leq}{\leqslant}

\newcommand{\KU}{\mathrm{KU}}

\DeclareMathOperator{\Br}{Br}

\newcommand{\Gm}{\mathds{G}_{m}}

\newcommand{\et}{\mathrm{\acute{e}t}}

\newcommand{\per}{\mathrm{per}}

\newcommand{\ind}{\mathrm{ind}}

\DeclareMathOperator{\Spec}{Spec}

\newcommand{\iso}{\cong}

\theoremstyle{plain}
\newtheorem{theorem}{Theorem}[section]
\newtheorem*{theorem*}{Theorem}

\newtheorem*{conjecture*}{Conjecture}

\newtheorem*{corollary*}{Corollary}

\theoremstyle{plain}

\theoremstyle{definition}

\newtheoremstyle{named}{}{}{\itshape}{}{\bfseries}{.}{.5em}{#1 \thmnote{#3}}
\theoremstyle{named}
\newtheorem*{namedtheorem}{Theorem}
\newtheorem*{namedconjecture}{Conjecture}
\newtheorem*{namedcorollary}{Corollary}

\theoremstyle{definition}

\newtheorem*{example*}{Example}

\newtheorem*{question*}{Question}

\begin{document}

\maketitle

\begin{abstract}
 \noindent
 We prove the topological analogue of the period--index conjecture in each
 dimension away from a small set of primes.\footnote{{\bf Key Words;} Brauer groups, period--index problems, twisted
 $K$-theory. {\bf Mathematics Subject Classification 2010:}
 \href{http://www.ams.org/mathscinet/msc/msc2010.html?t=14Fxx&btn=Current}{14F22},
 \href{http://www.ams.org/mathscinet/msc/msc2010.html?t=18Exx&btn=Current}{19L50},
 \href{http://www.ams.org/mathscinet/msc/msc2010.html?t=18Exx&btn=Current}{55R40}.}
\end{abstract}

\renewcommand{\abstractname}{R\'esum\'e}

\begin{abstract}
 \noindent
 Nous prouvons l'analogue topologique de la conjecture exposant et indice dans chaque dimension en dehors d'un petit ensemble de nombres premiers.
\end{abstract}

\section{Introduction}

We prove the following theorem, partly solving the so-called topological period--index conjecture of~\cite{aw3}. For
background, see~\S\ref{sec:background}.

\begin{namedtheorem}[A]\label{thm:a}
 Let $X$ be a finite $2d$-dimensional CW complex\footnote{More precisely, what we require is for $X$ to have the 
 homotopy type of a retract of a finite CW complex and for $\H^i(X,A)=0$ for $i>2d$ and any abelian group
 $A$.} and
 let $\alpha\in\Br(X)=\H^3(X,\ZZ)_{\tors}$ be a Brauer class. Setting $n=\per(\alpha)$, we have
 $$\ind(\alpha)|n^{d-1}\prod_{p|n}p^{v_p((d-1)!)},$$ where the product
 ranges over the prime divisors of $n$ and where $v_p$ denotes the $p$-adic valuation.
\end{namedtheorem}

Away from $(d-1)!$, the result simplifies.

\begin{namedcorollary}[B]\label{cor:b}
 If $X$ is a finite $2d$-dimensional CW complex, $\alpha\in\Br(X)$, and $\per(\alpha)$ is prime to $(d-1)!$, then
 $\ind(\alpha)|\per(\alpha)^{d-1}$.
\end{namedcorollary}

The corollary is an exact topological analogue, away from some small primes, of the well-known period--index conjecture
for division algebras over function fields (see~\cite[Section~2.4]{colliot}).

\begin{conjecture*}[Function field period--index conjecture]
 Let $k$ be algebraically closed and let $K$ be of transcendence degree $d$ over $k$. If $\alpha\in\Br(K)$,
 then $$\ind(\alpha)|\per(\alpha)^{d-1}.$$
\end{conjecture*}

Tsen's theorem implies that $\Br(K)=0$ when $K$ has
transcendence degree $1$ over an algebraically closed field
(see~\cite{gille-szamuely}). In~\cite{dejong}, de Jong proved
the conjecture when $d=2$ for Brauer classes of order
relatively prime to the characteristic. Lieblich removed this
restriction in~\cite{lieblich-twisted} thus establishing the $d=2$ case in full
generality. There are no other cases known of the conjecture. More precisely,
the period--index conjecture for function fields is not known for a single function field of
transcendence degree $d>2$ over an algebraically closed field.

Nevertheless, work of de Jong and Starr~\cite{dejong-starr} has reduced
the conjecture above to the following special cases.

\begin{conjecture*}[Global period--index conjecture]
 Let $k$ be an algebraically closed field and let $X$ be a smooth projective
 $k$-scheme of dimension $d$. If $\alpha\in\Br(X)\subseteq\Br(k(X))$, then
 $$\ind(\alpha)|\per(\alpha)^{d-1}.$$
\end{conjecture*}

In other words, to prove the period--index conjecture for function fields, it is
enough to prove it for unramified classes. We therefore view the period--index
conjecture as a global problem. If $k=\CC$, then the space of complex points $X(\CC)$ of
a smooth projective $\CC$-scheme of dimension $d$ admits the structure of a
$2d$-dimensional CW complex. Thus,
Corollary~\hyperref[cor:b]{B} provides evidence for this conjecture.

The topological period--index problem was introduced by the authors
in~\cite{aw1} where weak lower bounds were given and where the $d\leq 2$ cases
were solved. The $d=3$ case of the topological period--index problem was settled
in~\cite{aw3}, where it was proved moreover that the bound appearing in
Theorem~\hyperref[thm:a]{A} is sharp in that case. The case of $d=4$ was solved by Gu
in~\cite{gu2,gu3}, where the upper bound of Theorem~\hyperref[thm:a]{A} was
found independently (in the $d=4$ case) and where it is shown that the bound appearing in the theorem is
sharp for square-free classes. The best possible bound for $d=4$ and $n=\per(\alpha)$, as
proved by Gu, is
$$\ind(\alpha)|\begin{cases}
 e_3(n)n^3&\text{if $4|n$ and}\\
 e_2(n)e_3(n)n^3&\text{otherwise,}
\end{cases}$$
where $e_p(n)=p$ if $p|n$ and $1$ otherwise. In other words, there exist
$8$-dimensional finite CW complexes where these upper bounds on the index are
achieved.

The fact that for $d=3$ there exist $6$-dimensional finite CW complexes with
Brauer classes $\alpha$ having $\per(\alpha)=2$ and $\ind(\alpha)=8$ leads to
the natural question in~\cite{aw3} of whether the global period--index conjecture might be false
at the prime $2$ for threefolds over the complex numbers. However, recent work of Crowley--Grant~\cite{crowley-grant} proves
(a) that these topological examples can be found among closed orientable
$6$-manifolds but that (b) these examples cannot be found among closed $6$-dimensional
$\mathrm{Spin}^{\mathrm{c}}$-manifolds and hence they cannot be found among
closed $6$-manifolds of the form $X(\CC)$ for a smooth projective complex $3$-fold $X$.
The question of what happens for $d=4$ is the subject of ongoing work of
Crowley--Gu--Haesemeyer who prove in~\cite{crowley-gu-haesemeyer} that for closed orientable
$8$-manifolds $X$ one has $\ind(\alpha)|\per(\alpha)^3$ for $\alpha\in\Br(X)$ unless
$\per(\alpha)\equiv 2\pmod{4}$ in which case $\ind(\alpha)|2\per(\alpha)^3$.

Note that Gu's bound $\ind(\alpha)|e_3(n)n^3$ if $4|n$ is better than the bound
$\ind(\alpha)|e_2(n)e_3(n)n^3$ arising from Theorem~\hyperref[thm:a]{A}.
We do not further address in this paper the sharpness of the bounds in
Theorem~\hyperref[thm:a]{A} except to make the following conjecture.

\begin{namedconjecture}[C]\label{conj:c}
 The bounds of Corollary~\hyperref[cor:b]{B} are the best possible. That is,
 for every $d\geq 1$ and every natural number $n$ prime to $(d-1)!$, there exists
 a finite $2d$-dimensional CW complex $X$ and a Brauer class
 $\alpha\in\Br(X)$ such that $\per(\alpha)=n$ and $\ind(\alpha)=n^{d-1}$.
\end{namedconjecture}

The bounds in the period--index conjecture for function fields are known to be sharp,
for example by Gabber's appendix to~\cite{colliot-exposant}.

\paragraph{Acknowledgments.}
We would like to thank Diarmuid Crowley, Mark Grant, Xing Gu, and Christian
Haesemeyer for their
interest in this project and for many conversations over the years about
the topological period--index problem. The first author was supported by NSF Grant DMS-1552766.

\section{Background and strategy}\label{sec:background}

We quickly review the period--index problem in three settings.

\paragraph{Period--index for fields.}
The period--index problem originated in the domain of division algebras over fields. Specifically, for a field $K$, we
have the Brauer group $\Br(K)=\H^2_{\et}(\Spec K,\Gm)$. This group is isomorphic to the set of isomorphism classes of
finite dimensional division $K$-algebras with center exactly $K$. Given $\alpha\in\Br(K)$, we have two numbers:
$\per(\alpha)$, which is the order of $\alpha$ in the torsion abelian group $\Br(K)$, and $\ind(\alpha)$ which is the
unique positive integer such that $\ind(\alpha)^2=\dim_KD$ where $D$ is a division algebra with Brauer class
$[D]=\alpha$. It is not hard to see that
$$\per(\alpha)|\ind(\alpha)$$ and Noether proved that these two numbers have
the same prime divisors. It follows that there is some integer $e_\alpha$ such that
$\ind(\alpha)|\per(\alpha)^{e_\alpha}$.

The period--index problem is to find for a fixed field $K$ a number $e$
such that $$\ind(\alpha)|\per(\alpha)^e$$ for all $\alpha\in\Br(K)$ and, this
done, to find the smallest such number. For example, the
Albert--Brauer--Hasse--Noether theorem says that if $K$ is a number field, then
$\ind(\alpha)=\per(\alpha)$ so that $e=1$ works (see~\cite[Remark~6.5.6]{gille-szamuely}).

The period--index conjecture for function fields can be
rephrased as saying that if $K$ has transcendence degree $d$ over an
algebraically closed field, then $e=d-1$ is the solution. For $d\geq 3$, it is
not yet known that there is any $e$ that works for all Brauer classes, but
some results are known one prime at a time~\cite{matzri}.

\paragraph{Period--index for schemes.}
We introduced the period--index problem in other settings in~\cite{aw1}.
For instance, if $X$ is a quasicompact scheme, then the Brauer group
$$\Br(X)\subseteq\H^2_\et(X,\Gm)_\tors\subseteq\H^2_\et(X,\Gm)$$ of Azumaya
algebras of Grothendieck~\cite{grothendieck-brauer-1} is a torsion abelian group.
Given $\alpha\in\Br(X)$, we again let $\per(\alpha)$ be the order of $\alpha$
in $\Br(X)$. We define $$\ind(\alpha)=\gcd\{\deg(\Ascr)\colon\text{$\Ascr$ is
an Azumaya algebra with $[\Ascr]=\alpha$}\}.$$ In~\cite{aw4}, we showed that
even on smooth schemes over the complex numbers, it is necessary to take the
greatest common divisor to obtain a good theory. In this setting, we have
$\per(\alpha)|\ind(\alpha)$ and the numbers have the same prime divisors
by~\cite{aw7}. Thus, one can formulate the period--index problem for $X$.

\paragraph{Period--index for topological spaces.}
Finally, if $X$ is a topological space, we have the Brauer group
$\Br(X)\subseteq\H^3(X,\ZZ)_\tors\subseteq\H^3(X,\ZZ)$ of topological
Azumaya algebras, also introduced in~\cite{grothendieck-brauer-1}.
We define $\per(\alpha)$ and $\ind(\alpha)$ as for schemes. Again,
$\per(\alpha)|\ind(\alpha)$ and we proved
that these have the same prime divisors if $X$ is a finite CW complex
in~\cite{aw1} and in general in~\cite{aw7}.

The topological period--index problem of~\cite{aw1} asks the following. Given $d\geq 1$,
find bounds $e$ such that if $X$ is a finite $2d$-dimensional CW complex and
$\alpha\in\Br(X)$, then $\ind(\alpha)|\per(\alpha)^{e}$. We proposed $e=d-1$ as
a straw man in~\cite{aw3}, where we immediately proved that for $d=3$ this
fails in general when $2|\per(\alpha)$. Gu has proved this fails for $d=4$ if
$2$ or $3$ divides $\per(\alpha)$. But, in these low-dimensional cases, these
small primes are the only obstruction. We prove in Theorem~\hyperref[thm:a]{A}
that this pattern continues in higher dimensions.

The topological results reveal a pattern which has not yet been discovered in
algebra: a dependence on the prime divisors of the period and their
relationship to $d$. This dependence comes, as we will see, from the `jumps' in
the cohomology of the Eilenberg--MacLane spaces $K(\ZZ/(n),2)$.

\paragraph{Strategy.}
If $X$ is a topological space and $\alpha\in\H^3(X,\ZZ)$, Atiyah and Segal
constructed in~\cite{atiyah-segal-1} an
$\alpha$-twisted form of complex $K$-theory $\KU(X)_\alpha$ and in~\cite{atiyah-segal-2} an $\alpha$-twisted
Atiyah--Hirzebruch spectral sequence
$$\E_2^{s,t}=\H^s(X,\ZZ(\tfrac{t}{2}))\Rightarrow\KU^{s+t}(X)_\alpha,$$
where $\ZZ(\tfrac{t}{2})\iso\ZZ$ if $t$ is even and $\ZZ(\tfrac{t}{2})=0$ if $t$ is odd.
The differentials $d_r^\alpha$ have bidegree $(r,1-r)$.
We proved in~\cite{aw1} that if $X$ is a connected finite CW complex, then $\ind(\alpha)$
generates the group $$\E_\infty^{0,0}\subseteq\H^0(X,\ZZ)\iso\ZZ$$ of permanent
cycles. Thus, to bound the index, one attempts to bound the orders of the
differentials $$d_{2r+1}\colon\E_r^{0,0}\rightarrow\E_r^{2r+1,-2r}.$$
There is a universal case to consider for all order $m$ topological Brauer classes, namely
the space $K(\ZZ/(m),2)$ and a generator $\alpha$ of $\H^3(K(\ZZ/(m),2),\ZZ)\iso\ZZ/(m)$.
By studying the orders of the differentials in this particular case,
we prove Theorem~\hyperref[thm:a]{A}.

% \begin{remark}[$2d$-dimensional versus $(2d-1)$-dimensional]
% The connection between the topological period--index problem and the differentials in the
% twisted Atiyah--Hirzebruch spectral sequence implies that the general solution to the
% topological period--index
% problem is the same for finite $2d$-dimensional CW-complexes as for finite
% $(2d-1)$-dimensional CW complexes.
% \end{remark}

\section{The cohomology of $K(\ZZ/(n),2)$}\label{sec:cartan}

We recall some results of Cartan~\cite{cartan} on the cohomology of Eilenberg--MacLane spaces in the special case of
$K(\ZZ/(n),2)$, which we will use in the next section to give upper bounds on the index of period $n$ classes. We claim no
originality in our presentation here, but we hope its inclusion will be useful to the reader.

Write
\begin{equation*}
 n=p_1^{r_1}\cdots p_k^{r_k},
\end{equation*}
and pick a generator $u_i$ of the subgroup $\ZZ/(p_i^{r_i})$ of $\ZZ/(n)$ for $1\leq i\leq k$.
Let $v_i=p_i^{r_k-1}u_i$.
Cartan~\cite[Th\'eor\`eme 1]{cartan} gives a recipe for computing the integral homology (Pontryagin) ring of
$K(\ZZ/(n),2)$. 

For each prime $p$ and positive integer $f$, consider certain words in the $4$ symbols:
\begin{equation*} \sigma, \quad \gamma_p, \quad \phi_p, \quad \psi_{p^f}\end{equation*}
The symbol $\psi_{p^f}$, if it appears, is the last symbol in a word. The \textit{height} of a word $\alpha$ is the
total number of $\sigma$, $\phi_p$ and $\psi_{p^f}$ appearing. The \textit{degree} of
$\alpha$ is defined recursively by letting $\deg(\emptyset)=0$ and
\begin{equation*}
 \deg(\sigma\alpha)=1+\deg(\alpha),\hspace{1cm}\deg(\gamma_p\alpha)=p\deg(\alpha),\hspace{1cm}\deg(\phi_p\alpha)=2+p\deg(\alpha),
 \hspace{1cm} \deg(\psi_{p^f})= 2
\end{equation*}
For each prime $p$, an \textit{admissible} $p$-word $\alpha$ is a word
on $3$ symbols $\sigma$, $\gamma_p$, and $\phi_p$ such that $\alpha$ is non-empty, the
first and last letters of $\alpha$ are $\sigma$ or $\phi_p$, and for each letter
$\gamma_p$ or $\phi_p$, the number of letters $\sigma$ appearing to the right in $\alpha$ is
even. % The \textit{height} of $\alpha$ is the total number of $\sigma$ and $\phi_p$. The \textit{degree} of
% $\alpha$ is defined recursively by
% \begin{equation*}
% \deg(\sigma\alpha)=1+\deg(\alpha),\hspace{1cm}\deg(\gamma_p\alpha)=p\deg(\alpha),\hspace{1cm}\deg(\phi_p\alpha)=2+p\deg(\alpha)
% \end{equation*}
% starting with $\deg(\sigma)=1$ and $\deg(\phi_p)=2$.
% % Admissible $p$-words are divided into those of the first type, where the last letter is
% % $\sigma$, and those of the second type, where the last letter is $\phi_p$.
In addition to the admissible words, we will use the auxiliary words $\sigma^{h-1}\psi_{p^f}$, of height $h$ and
degree $h+1$.

Let $E(x,2q-1)$ denote the exterior graded algebra over $\ZZ$ with generator $x$ of degree
$2q-1$ endowed with the trivial dg-algebra structure $dx=0$. Let $P(x,2q)$ be the divided power polynomial
algebra over $\ZZ$ with generator $x$ of degree $2q$ and given the trivial dg algebra
structure with $dx=0$. Cartan calls $E(x,2q-1)$ and $P(x,2q)$
\textit{elementary complexes of the first type}. Define tensor dg algebras
$E(x,2q-1)\otimes_\ZZ P(y,2q)$ by $dx=0$ and $dy=hx$
for some integer $h$ and $P(x,2q)\otimes_\ZZ E(y,2q+1)$ by $dx=0$ and $dy=hx$ (the integer $h$
 is part of the data, even though it is not specified in the notation).
These are the \textit{elementary complexes of the second type}.
The positive-degree homology groups of $E(x,2q-1)\otimes_\ZZ P(x,2q)$ are
\begin{equation*}
 \H_{2q-1+2qk}(E(x,2q-1)\otimes_\ZZ P(x,2q))=\ZZ/h\cdot x\gamma_k(y)
\end{equation*}
for $k\geq 0$ and
and $0$ otherwise, where $\gamma_k$ is the $k$th divided power operation. For
$P(x,2q)\otimes_\ZZ E(y,2q+1)$, we get
\begin{equation*}
 \H_{2qk}(P(x,2q)\otimes_\ZZ E(y,2q+1))=\ZZ/hk\cdot\gamma_k(x)
\end{equation*}
for $k\geq 0$
and all other homology groups are $0$.

The height-$2$ admissible or auxiliary $p$-words are
\begin{equation*}
 \sigma^2,\hspace{1cm}\sigma\gamma^k_p\phi_p,\hspace{1cm}\phi_p\gamma^k_p\phi_p,\hspace{1cm}\sigma\psi_{p^f}
\end{equation*}
for $k\geq 0$ and $f\geq 1$. These are of degrees $2$, $1+2p^k$, $2+2p^{k+1}$, and $3$, respectively. Below, the symbols
$ u_i$ and $v_i$, are just formal indeterminates to keep track of generators for different dg algebras.

For each $p_i$, we define a dg algebra $X_{p_i}$ as the following tensor product of elementary complexes of
the second type
\begin{gather*}
 X_{p_i}=P(\sigma^2u_i,2)\otimes E(\sigma\psi_{p_i^{r_i}} u_i,3)\\
% \bigotimes_{k=0}^{\infty}E(\sigma\gamma_{p_i}^{k+1}\phi_{p_i} u_i,1+2{p_i}^{k+1})\otimes
% P(\phi_{p_i}\gamma_{p_i}^k\phi_{p_i} u_i,2+2{p_i}^{k+1})\\
 \bigotimes_{k=0}^{\infty}E(\sigma\gamma_{p_i}^{k+1}\phi_{p_i} v_i,1+2{p_i}^{k+1})\otimes
 P(\phi_{p_i}\gamma_{p_i}^k\phi_{p_i} v_i,2+2{p_i}^{k+1}).
\end{gather*}
The differentials are
\begin{equation*}
 d(\sigma\psi_{p_i^{r_i}} u_i)=p_i^{r_i}\sigma^2u_i\quad\text{and}\quad
% \\\hspace{9mm} d(\phi_{p_i}\gamma_{p_i}^k\phi_{p_i} u_i)=p_i^{r_i}\sigma\gamma_{p_i}^{k+1}\phi_{p_i} u_i,\hspace{9mm}
 d(\phi_{p_i}\gamma_{p_i}^k\phi_{p_i} v_i)=p_i\sigma\gamma_{p_i}^{k+1}\phi_{p_i} v_i,
\end{equation*}
i.e., $h=p_i^{r_i}$
% \\, $h=p_i^{r_i}$,\\
or $h=p_i$, respectively. (Here the $r_i$ are the
exponents appearing in the prime decomposition of $n$.)
Let
\begin{equation*}
 X=X_{p_1}\otimes\cdots\otimes X_{p_k}.
\end{equation*}
Then, Cartan~\cite{cartan}*{Th\'eor\`eme 1} gives a surjection, which depends on the choice
of the $u_i$,
\begin{equation*}
 \H_k(X)\rightarrow\H_k(K(\ZZ/(n),2),\ZZ),
\end{equation*}
and the kernel is described. This map induces for each $i$ a surjection
\begin{equation*}
 \H_k(X_{p_i})\rightarrow\H_k(K(\ZZ/(n),2),\ZZ)\{p_i\},
\end{equation*}
the $p_i$-primary part of homology. The largest possible torsion in $\H_{2k}(X_{p_i})$ comes from the
first term $P(\sigma^2 u_2,2)\otimes E(\sigma\psi_{p_i^{r_i}} u_i,3)$. Specifically, by using the
K\"unneth theorem, the exponent of $\H_{2k}(X_{p_i})$ is the same as the exponent of
$\H_{2k}(P(\sigma^2 u_2,2)\otimes E(\sigma\psi_{p_i^{r_i}} u_i,3))$, namely
\begin{equation*}
 p_i^{r_i}k.
\end{equation*}
It follows that the exponent of $\H_{2k}(X)$ is at most $nk$. Write $k=k'a$,
where $a$ is largest integer dividing $k$ and prime to $m$. Since
$\H_{2k}(K(\ZZ/(n),2),\ZZ)$ is $n$-primary torsion, we see that $nk'$ kills all of
$\H_{2k}(K(\ZZ/(n),2),\ZZ)$; in fact, the exponent is exactly $nk'$, as follows from Cartan's
description of the kernel of $\H_{2k}(X)\rightarrow\H_{2k}(K(\ZZ/(n),2),\ZZ)$, but we will not
need this fact.

Since we will be able to argue prime-by-prime in a moment, it is helpful to record the $p_i$-primary part of the exponent
of $\H_{2k}(K(\ZZ/(n), 2), \ZZ)$ for a given prime $p_i$. If $\xi \in \H_{2k}(K(\ZZ/(n), 2), \ZZ)$ is an element of
$p_i$-primary order, then the order of $\xi$ divides $p_i^{r_i + v_{p_i}(k)}$.

\section{Proof of the main theorem}

If $X$ is a topological space with finitely many connected components, then the topological Brauer group $\Br(X)$ is a
subgroup of $\Br'(X)=\H^3(X,\ZZ)_{\tors}$. Serre showed that if $X$ is compact, then $\Br(X)=\Br'(X)$
(see~\cite{grothendieck-brauer-1}). This will be the main setting of this paper. We want to make a general remark on a
different version of the index to which our methods apply whether or not $X$ is compact.

We may weaken the
hypothesis on the finiteness of $X$ in Theorem~\hyperref[thm:a]{A} at the expense of using the $K$-theoretic index rather than the
topological index. Recall from~\cite{aw1} that the $K$-theoretic index $\ind_\K(\alpha)$ is
defined as the (positive) generator of the image of the rank map
$\KU^0(X)_\alpha\rightarrow\ZZ$, where $\KU^0(X)_\alpha$ denotes the
$\alpha$-twisted $K$-theory group. When $X$ is finite-dimensional (and
connected), one computes
$\ind_\K(\alpha)$ as the generator of the group
$$\E_\infty^{0,0}\subseteq\E_2^{0,0}\iso\ZZ$$
by the convergence of the $\alpha$-twisted Atiyah--Hirzebruch spectral
sequence;
hence the differentials
$d_{2k+1}\colon\E_{2k+1}^{0,0}\rightarrow\E_{2k+1}^{2k+1,-2k}$ of the twisted
Atiyah--Hirzebruch spectral sequence control the
$K$-theoretic index. It is this crucial fact we exploit below.

In general, $\per(\alpha)|\ind_\K(\alpha)|\ind(\alpha)$ and if $X$ is compact,
then $\ind_\K(\alpha)=\ind(\alpha)$. Thus, Theorem~\hyperref[thm:a]{A} follows
from the following theorem.

\begin{theorem}
 Let $X$ be a $2d$-dimensional CW complex, and let $\alpha\in\Br'(X)$ have period
 $m=p_1^{r_1}\cdots p_k^{r_k}$ where the $p_i$ are distinct primes. Then,
 \begin{equation*}
 \ind_\K(\alpha)|\prod_{i=1}^k
 {p_i}^{(d-1)r_i+v_{p_i}(2)+\cdots+v_{p_i}(d-1)}=m^{d-1}\prod_{i=1}^k {p_i}^{v_{p_i}((d-1)!)},
 \end{equation*}
 where $v_{p_i}$ is the $p_i$-adic valuation.
 \begin{proof}
According to the main result of~\cite{aw8}, it is enough to prove the theorem when $m$ is a prime power. So, suppose that
 $m=p_1^{r_1}=p^r$. Let $\beta$ be a generator of
 $\H^3(K(\ZZ/(p^r),2),\ZZ)\iso\ZZ/(p^r)$. There is some map $\sigma\colon
 X\rightarrow K(\ZZ/(p^r),2)$ such that $\sigma^*\beta=\alpha$. The twisted Atiyah--Hirzebruch spectral sequence is
 functorial, so we obtain---in particular---a map on the $(0,0)$-terms in each page:
 \[ \E^{0,0}_j (K(\ZZ/(p^r),2)) \to \E^{0,0}_j(X).\]
An elementary induction
argument shows that this map is a monomorphism on each page, and so by cellular approximation, $\ind_\K(\alpha)$ is bounded above by $\ind_\K$
 for the restriction of $\beta$ to a $2d$-skeleton of $K(\ZZ/(n), 2)$. 
% if
% $n_j$ generates $\E_{2j+1}^{0,0}(K(\ZZ/(p^r),2),\beta)\subseteq\ZZ$ in the $\beta$-twisted
% Atiyah--Hirzebruch spectral sequence for $K(\ZZ/(p^r),2)$, then $n_j$ is
% also in $\E_{2j+1}^{0,0}(X,\alpha)\subseteq X$ and if
% $Nd_{2j+1}(n_j)=0$ in $\E_{2j+1}^{2j+1,-2j}(K(\ZZ/(p^r),2),\beta)$, then
% $Nd_{2j+1}(n_j)=0$ in $\E_{2j+1}^{2j+1,-2j}(X,\alpha)$. Thus, to prove
% the theorem, it is enough to consider the $2d$-skeleton of
% $K(\ZZ/(p^r),2)$.

 There is an isomorphism
 \begin{equation*}
 \tilde \H^{2j+1}(K(\ZZ/(p^r),2),\ZZ)\iso \tilde\H_{2j}(K(\ZZ/(p^r),2),\ZZ).
 \end{equation*}
 We have seen in Section~\ref{sec:cartan} that the latter group is $p^r p^{v_p(j)}$-torsion. In the Atiyah--Hirzebruch
 spectral sequence for the Bockstein $\beta$ on a $2d$-skeleton of $K(\ZZ/(p^r),2)$, only the differentials
 $d_{2j+1}^{\beta}$ for $1\leq j\leq d-1$ are possibly non-zero---in particular, the cohomology group in degree $2d$ of
 the skeleton, which may differ from that of $K(\ZZ/(p^r), 2)$, plays no part in the calculation. We are interested in
 the order of the image of
 \begin{equation*}
 d_{2j+1}^{\beta}\colon:\E_{2j+1}^{0,0}\rightarrow\E_{2j+1}^{2j+1,-2j},
 \end{equation*}
 where the latter group is a subquotient of $\H^{2j+1}(K(\ZZ/(p^r),2),\ZZ)$, and
 hence the order of the image divides $p^{r+v_p(j)}$. This gives the result.
 \end{proof}
\end{theorem}

\small
\bibliographystyle{amsplain}
\bibliography{tpic}

\printindex

% \vspace{20pt}
% \scriptsize
% \noindent
% Benjamin Antieau\\
% University of Illinois at Chicago\\
% Department of Mathematics, Statistics, and Computer Science\\
% 851 South Morgan Street, Chicago, IL 60607\\
% \texttt{benjamin.antieau@gmail.com}

\end{document}